\documentclass[leqno, a4paper, 11pt]{article}

\usepackage{amsfonts, amssymb, calrsfs}
\usepackage{latexsym}

\def\rk{\mathop{\mathrm{rk}}\nolimits}
\def\mod{\mathop{\mathrm{mod}}\nolimits}
\def\num{\mathop{\mathrm{wn}}\nolimits}
\def\den{\mathop{\mathrm{wd}}\nolimits}

\def\Z{\mathbf{Z}}
\def\Q{\mathbf{Q}}
\def\C{\mathbf{C}}

\def\O{\mathcal{O}}
\def\G{\mathcal{G}}
\def\p{{\mathfrak{p}}}
\def\q{{\mathfrak{q}}}
\def\ker{\mathop{\mathrm{ker}}\nolimits}
\def\til#1{\widetilde{#1}}

\def\pf{{\indent\textit{Proof.}\ }}
\def\qed{\hfill$\square$}
\newcounter{para}[section]
\renewcommand{\thepara}{\thesection.\arabic{para}}
\renewcommand{\thesection}{\arabic{section}}
\renewcommand{\paragraph}{\refstepcounter{para}
\indent{\bf{\thepara}}}
\def\sectioning#1{\vv 

 \refstepcounter{section}
\indent {\bf \thesection. #1}}
\newcommand{\vv}{\vspace{1ex}}

\begin{document}

{\footnotesize \noindent  Running head: HTP for rings of integers \hfill 
December 19th, 2003 

\noindent Math.\ Subj.\ Class.\ (2000): 03B25, 11U05  }

\vspace{4cm}

\begin{center} 

{\Large Division-ample sets and 

\vv

the Diophantine problem for rings of integers}

\vv

{\sl by} Gunther Cornelissen, Thanases Pheidas {\sl and} Karim Zahidi

\vv

\end{center}

\vv

\vv

\begin{quote}
{\small {\bf Abstract.} 
We prove that Hilbert's Tenth Problem for a ring of integers in a number field 
$K$ has a negative answer if $K$ satisfies two arithmetical conditions 
(existence of a so-called {\sl division-ample} set of integers and of an elliptic 
curve of rank one over $K$). We relate division-ample sets to arithmetic 
of abelian varieties.}

\end{quote}

\vv

\vv

{\bf Introduction.} 

\vv

Let $K$ be a number field and let $\O_K$ be its ring of integers. {\sl Hilbert's 
Tenth Problem} or {\sl the diophantine problem} for $\O_K$ is the following: is 
there an algorithm (on a Turing machine) that decides whether an arbitary 
diophantine equation with coefficients in $\O_K$ has a solution in $\O_K$. 

The answer to this problem is known to be negative if $K=\Q$ (\cite{Davis}), 
and for several other such $K$ (such as imaginary quadratic number fields 
\cite{Denef}, totally real fields \cite{Denef2}, abelian number fields 
\cite{ShapiroShlap}) by reduction to the case $K=\Q$. This reduction consists in 
finding a {\sl diophantine model} (cf.\ \cite{CornelissenZahidi}) for integer 
arithmetic over $\O_K$. The problem is open for general number fields 
(for a survey see \cite{PhZ} and \cite{Shl}), but has been solved conditionally, e.g.\ by Poonen \cite{Poonen} (who shows that 
the set if rational
integers is diophantine over $\O_K$ if there exists an elliptic curve over $\Q$ 
that has rank one over both $\Q$ and $K$). In this paper, we give a more general 
condition as follows:

\vv

{\bf Theorem.} \ {\sl The diophantine problem for the ring of integers $\O_K$ of 
a number field $K$ has a negative answer if the following exist:

(i) an elliptic curve defined over $K$ of rank one over $K$;

(ii) a \emph{division-ample} set $A \subseteq \O_K$.}

\vv

\noindent A set $A \subseteq \O_K$ is called {\sl division-ample} if the 
following three conditions are satisfied:

$\diamond$ {\sf (diophantineness)} $A$ is a diophantine subset of $\O_K$;

$\diamond$ {\sf (divisibility-density)} Any $x \in \O_K$ divides an element of 
$A$; 

$\diamond$ {\sf (norm-boundedness)} There exists an integer $\ell>0$, such that 
for any $a \in A$, there is an integer $\til{a} \in \Z$ with 
$\til{a}$ dividing $a$ and $|N(a)| \leq |\til{a}|^\ell$.

\vv

{\bf Proposition.} \ {\sl A division ample set exists if either 

(i) there exists an abelian 
variety $G$ over $\Q$ such that $$ \rk G(\Q) = \rk G(K) > 0; \mbox{ or } $$

(ii) there exists a 
commutative (not necessarily complete) 
group variety $G$ over $\Z$ such that $G(\O_K)$ is finitely generated and such 
that
$ \rk G(\Q) = \rk G(K) > 0. $}

\vv

From (i) in this proposition, it follows that our theorem includes that of Poonen, but it isolates
the notion of ``division-ampleness'' and shows it can be satisfied in a broader context. It would for example be interesting to construct, for a given number field $K$, a curve over $\Q$ such that it's Jacobian satisfies this condition.

As we will prove below, part (ii) of this proposition is satisfied for the relative norm one torus
$G=\ker(N_K^{KL})$ for a number field $L$ linearly disjoint from $K$, if $K$ is quadratic imaginary 
(choosing $L$ totally real). 

It  would be interesting to know other division-ample sets, in particular, 
such that are not subsets of the integers.

The proof of the main theorem will use divisibility on elliptic curves and a 
lemma
from algebraic number theory of Denef and Lipshitz. Some of our arguments are similar to ones in \cite{Poonen}, but 
we have avoided continuous reference both for reasons of completeness 
and because our results have been obtained independently.

\vv

\vv

\sectioning{Lemmas on number fields} \label{nf}

\vv

In this section we collect a few facts about general number fields which 
will play a r\^ole in subsequent proofs. Fix $K$ to be a number field, let 
$\O=\O_K$ be its ring of integers, and let $h$ denote the class number of $\O$.
Let $N=N^K_\Q$ be the norm from $K$ to $\Q$, and let $n=[K:\Q]$ denote the 
degree of $K$. Let $|$ denote ``divides'' in $\O$. 

First of all, we will say a subset $S \subseteq K^n$ is ``diophantine over 
$\O$'' if its set of representatives $ \til{S} \subseteq (\O \times 
(\O-\{0\}))^n$ given by
$$ \til{S} := \{ (a_i,b_i)_{i=1}^n \in (\O \times (\O-\{0\}))^n \ | \ 
({a_i}/{b_i})_{i=1}^n \in S \}$$
is diophantine over $\O$. Recall that ``$x \neq 0$'' is diophantine over $\O$ 
(\cite{DenefLipshitz} Prop.\ 1(b)), hence $S$ is diophantine over $\O$ if and 
only if it is diophantine over $K$.

Recall that there is no unique factorisation in general number fields, but we 
can use the following
valuation-theoretic remedy:

\vv
  
\paragraph\label{nf-1} {\bf Definition.} \ Let $x \in K$. If $x^h = \frac{a}{b}$ 
for $a,b \in \O$ with $(a,b)=1$ (the ideal generated by $a$ and $b$), we say that $a=\num(x)$ is {\sl a weak 
numerator} and $b=\den(x)$ is {\sl a weak denominator} for $x$.

\vv

\paragraph\label{nf-2} {\bf Lemma.} \ {\sl (i) For any $x \in K$ a weak 
numerator and 
a weak denominator exists and is unique up to units.

(ii) for any valuation, $v(x)>0 \iff v(\num(x))>0$ [and then $v(\num(x)) = h 
v(x)$], and $v(x)<0 \iff v(\den(x))>0$ [and then $v(\den(x))=-h v(x)$].

(iii) For $a \in \O, x \in K$, ``$a=\num(x)$'' and ``$a=\den(x)$'' are diophantine 
over $\O$.} 

\vv

\pf Since $\O$ is a Dedekind ring, $(x)$ has a unique factorisation in 
fractional ideals $$(x)=\p_1 \cdots \p_r \cdot \q^{-1}_1 \cdots \q^{-1}_s.$$ We 
let
$a$ be a generator for the principal ideal $(\p_1 \cdots \p_r)^h$ and $b$ a 
generator for $(\q_1 \cdots \q_s)^h$; these are obviously weak 
numerator/denominator for $x$. Uniqueness, (ii) and (iii) are obvious. \qed

\vv

\paragraph\label{nf-3} {\bf Lemma.} (Denef-Lipshitz \cite{DenefLipshitz}) {\sl  

(i) If $u \in \Z-\{0\}$ and $\xi \in \O$ satisfy the divisibility condition
$$ 2^{n!+1} \prod\limits_{i=0}^{n!-1} (\xi+i)^{n!}  \ | \ u $$
then for any embedding $\sigma : K \hookrightarrow \C$
$$ (\ast)_u \ \ \ \ \  |\sigma(\xi)| \leq \frac{1}{2} \sqrt[n!]{|N(u)|}. $$

(ii) If $\til{u} \in \Z-\{0\}, q \in \Z$ and $\xi \in \O$ satisfy 
$(\ast)_{\til{u}}$
for any embedding $\sigma : K \hookrightarrow \C$ and 
$ \xi \equiv q \mod \til{u}, $
then $\xi \in \Z$. }

\vv

\pf Easy to extract from the proof of Lemma 1 in \cite{DenefLipshitz}. \qed

\vv

\vv

\sectioning{Lemmas on elliptic curves} \label{ec}

\vv

Let $E$ denote an elliptic curve of rank one over $K$, written in Weierstrass 
form as 
$$ E \ : \ y^2 + a_1 xy + a_3 y = x^3 + a_2x^2 + a_4 x + a_6, $$ let $T$ be the 
order of the torsion group of $E(K)$, and let $P$ be a generator for the free 
part of $E(K)$. Define $x_n,y_n \in K$ by $nP=(x_n,y_n)$.

\vv

\paragraph\label{ec-1} {\bf Lemma.} \  {\sl For any integer $r$ the set
$rE(K)$ is diophantine over $K$ and, if $r$ is divisible by $T$, then $rE(K) = 
\langle rP \rangle \cong (\Z,+)$.} 

\vv

\pf A point $(x,y) \in K \times K$ belongs to $rE(K)-\{0\}$ if and only if 
$ \exists \ (x_0,y_0) \in E(K) \ : \ (x,y) = r(x_0,y_0)$. As the addition 
formul{\ae} on $E$ are algebraic with coefficients from $K$, this is a 
diophantine relation. The last statement is obvious. \qed

\vv

\paragraph\label{ec-2} {\bf Lemma.} \  {\sl There exists an integer $r>0$ such 
that for any non-zero integers $m,n \in \Z$, $m$ divides $n$ if and only if 
$\den(x_{rm}) | \den(x_{rn})$.} 

\vv

\pf We reduce the claim to a statement about valuations using lemma 
\ref{nf-2}(ii). The theory of the formal group associated to $E$ implies that if 
$n=mt$ and $v$ is a finite valuation of $K$ such that $v(x_{rm})<0$, then
$v(x_{rmt}) = v(x_{rm}) - 2 v(t) \leq v(x_{rm})$ (\cite{Silverman} VII.2.2).  

For the converse, we start by choosing $r_0$ in such a way that $r_0P$ is 
non-singular modulo all valuations $v$ on $K$. By the theorem of Kodaira-N\'eron 
(\cite{Silverman}, VII.6.1), such $r_0$ exists and it actually suffices to take 
$r_0=4 \prod_v v(\Delta_E)$, where $\Delta_E$ is the minimal discriminant of 
$E$, and the product runs over all finite valuations on $K$ for which 
$v(\Delta_E) \neq 0$. 
Note that then, $v(x_{r_0n})<0 \iff r_0nP=0$ in the group $E_v$ of non-singular 
points
of $E$ modulo $v$. 

We claim that for an arbitrary finite valuation $v$ on $K$, if $v(x_{r_0 n})<0$ 
and $ v(x_{r_0 m})<0$, then $v(x_{(r_0m,r_0n)})<0$, where $(\cdot,\cdot)$ 
denotes the gcd in $\Z$. Indeed, the hypothesis means $r_0mP=r_0nP=0$ in $E_v$. 
Since there exist
integers $a,b \in \Z$ with $(r_0m,r_0n)=ar_0m+br_0n$, we find $(r_0m,r_0n)P=0$ 
in $E_v$, and hence the claim. 

The main theorem of \cite{CheonHahn} states that for any sufficiently large $M 
(\geq M_0)$, there exists a finite valuation $v$ such that
$v(x_{M})<0$ but $v(x_i) \geq 0$ for all $i<M$. We choose $r=r_0 M_0$. Pick 
such a valuation $v$ for $M=rm$. The hypothesis implies that $v(x_{rn})<0$ and 
hence
$v(x_{r(m,n)})<0$. But $r(m,n) \leq rm$ 
and $v(x_{i})\geq 0$ for any $i<rm$. Hence $r(m,n)=rm$ so $m$ divides $n$. \qed

\vv

\paragraph\label{ec-3} {\bf Lemma.} \  {\sl Any $\xi \in \O-\{0\}$ divides the 
weak
denominator of some $x_n$.}

\vv

\pf Let $v$ be a valuation lying over the prime $p_v$ of $\Z$, and assume that 
$v(\xi)=e_v>0$. Then the group of non-singular points on $E$ modulo $v$ is 
finite, hence there exists an $n_v$ such that 
$n_vrP = 0$ in this group, i.e., $v(x_{n_vr})<0$. By the formal group law
formula, we see that $v(x_{p_v^{e_v}n_vr}) = v(x_{n_vr}) - 2 v(p_v^{e_v}) < 
-e_v$. Putting $$n=r \cdot \prod_{v(\xi)>0} n_v p_v^{e_v}$$
will then certainly suffice. \qed

\vv


\vv

\paragraph\label{ec-4} {\bf Lemma.} \  {\sl Let $m,n,q$ be integers with $n=mq$. 
Then $$ \den(x_m) | \num(\frac{x_n y_m}{y_n x_m}-q).$$}


\pf The formal power series expansion for addition on $E$ around $0$ 
(\cite{Silverman}, IV.2.3) implies that 
$ \frac{x_n}{y_n} = q \frac{x_m}{y_m} + O((\frac{x_m}{y_m})^2)$, from which the 
result follows. \qed

\vv

\vv

\sectioning{Proof of the main theorem} \label{mt}

\vv

Let $\xi \in \O$. Given an elliptic curve $E$ of rank one over $K$ as in the 
main theorem, we use the notation from section \ref{ec} for this $E$ --- in 
particular, choose a suitable $r$ such that lemma \ref{ec-2} applies; we also 
choose $\ell$ which comes with the definition of $A$. We claim that the 
following formul{\ae}  give a diophantine
definition of $\Z$ in $\O$:

$$\xi \in \Z \ \iff \ \exists \ m,n \in rT\Z, \ \exists \ u \in A - \{ 0 \} \ \ 
\left\{ \begin{array}{cl}
(1) & m|n \\

(2) & 2^{n!+1} \prod\limits_{i=0}^{n!-1} (\xi^{\ell n!}+i)^{n!} | u \\

(3) & u^h | \den(x_m) \\

(4) &  \den(x_m) | \num(\frac{x_n y_m}{y_m x_n}-\xi) 
\end{array} \right.$$

\paragraph\label{pf1} {\bf Any $\xi \in \Z$ satisfies the relations.} \ If $\xi \in 
\Z$, then a $u$ satisfying (2) exists because $A$ is division-dense. By lemma 
\ref{ec-3}, there exists an $m$ satisfying (3) for this $u$. Define $n=m \xi$ 
for this $m$. Then (1) is automatic and (4) is the contents of lemma \ref{ec-4}. 

\vv

\paragraph\label{pf2} {\bf A $\xi$ satisfying the relations is integral.} \ Let $q 
\in \Z$ satisfy $n=qm$ (which exists by (1)). Then lemma \ref{ec-4} implies that
$$ \den(x_m) | \num(\frac{x_n y_m}{x_m y_n}-q),$$
which can be combined with (4) using the non-archimedean triangle inequality to 
give
$$ \den(x_m) | \num(\xi-q) = (\xi-q)^h. $$
By (3), then also $u|\xi-q$. 

By norm-boundedness of $A$ we can find $\til{u} \in \Z$ such that $\til{u} | u$ 
and $|N(u)| \leq \til{u}^\ell.$ We still have
$$ (*) \ \ \xi \equiv q \mod \til{u}; \ \ \til{u},q \in \Z. $$
Condition (2) implies that Lemma \ref{nf-3}(i) can be applied with $\xi^{\ell 
n!}$ in place of $\xi$, so for any complex embedding $\sigma$ of $K$ we find
$$ (**) \ \ |\sigma(\xi)| \leq \frac{1}{2} | N(u)|^{\frac{1}{\ell n!}} \leq 
\frac{1}{2} N(\til{u})^{\frac{1}{n!}}. $$

Because of $(*)$ and $(**)$, we can apply Lemma \ref{nf-3}(ii) to conclude $\xi 
\in \Z$. 

\vv

\paragraph\label{pf3} {\bf The relations (1)-(4) are diophantine over $\O$.} \
By \ref{nf-2} and \ref{ec-1}, for $a \in \O$, the relations $\exists n \in rT\Z 
: a = \num(x_n)$ and $\exists n \in rT\Z : a = \den(x_n)$ are diophantine. By 
the diophantineness of $A$, the membership $u \in A$ is diophantine, and $u \neq 
0$ is diophantine (\cite{DenefLipshitz}, Prop.\ 1(b)). Condition (1) is 
diophantine because of Lemma \ref{ec-2}. Conditions (2)-(4) are obviously 
diophantine using \ref{nf-2}. \qed

\vv

\vv

\sectioning{Proof of the proposition and discussion of division-ample sets} \label{da}

\vv

\paragraph\label{da1} {\bf Rank-preservation over $\Q$.} We use \cite{CS} as a general reference on abelian varieties and formal groups. Suppose there exists an 
abelian variety $G$ of dimension $d$ over $\Q$ such that
$ \rk G(\Q) = \rk G(K) > 0$ (note that $G(K)$ is finitely generated by the Mordell-Weil theorem).  Let $T$ denote the (finite) order of the torsion of $G(K)$
and consider the free group $TG(K) \cong \Z^r$. The assumption implies that 
$G(\Q)$ is of finite index $[G(K):G(\Q)]$ in $G(K)$. The choice of an ample line 
bundle on $G$ gives rise to a projective embedding of $G$ in some projective
space with coordinates $\langle x_i \rangle_{i=1}^N$, where $G$ is cut out by 
finitely many polynomial equations and the addition on $G$ is algebraic in those 
coordinates. Suppose $\{t_i\}$ are algebraic function of the coordinates, and 
local uniformizers at the unit ${\bf 0} = (1:0:\dots:0)$ of $G$ (i.e., 
$\hat{\O}_{{G},{{\bf 0}}} = \Q[[t_1,\dots,t_d]]$), and define
$$ A_G := \{ \den(\prod_{i=2}^N t_i(P)) \ : \ P \in T [G(K):G(\Q)] \cdot G(K) \mbox{ and } t_1(P)=1 \}. $$
We claim that $A_G$ is division-ample. Indeed, the three conditions are satisfied:

\vv

(a) $A_G$ is obviously diophantine 
over $\O$ (the diophantine definition comes from the chosen embedding of $G$). 

\vv

(b) The analogue of lemma \ref{ec-3} remains valid: 

\begin{quote}
{\bf Claim.} \ {\sl $A_G$ is divisibility-dense.}
\end{quote}

\pf Since any $\xi \in \O - \{ 0 \}$ divides its norm, it suffice to prove
that any integer $z \in \Z - \{ 0 \}$ divides an element of $A_G$. Given 
a minimal model $\G_p$ for $G$ over a $p$-adic field $\Q_p$, let $\G_{p,0}$ 
denote the group of
points whose reduction is non-singular modulo $p$. Then $\G_{p,0}$ is a clopen 
subset of $\G_p$, so $\G_p/\G_{p,0}$ is finite (and non-trivial only for the 
finite set of primes for which $\G_p$ has bad reduction). Hence we can choose a 
finite $r$ so large that $rP_i$ is non-singular modulo all primes for all 
generators $P_i$ of $G(\Q)$. Pick a prime $p|z$, then since the group of
non-singular points on $G$ modulo $p$ is finite, there exists $n_p$ such that
$n_p rP = 0$ in this group, i.e., $v_p(t_i(P))>0$ for all $i$. The formal
group $\hat{G}_{{\bf 0}}$ of $G$ around ${\bf 0}$ (defined by the power series 
that give the addition in terms of $\{ t_i \})$ is a formal torus in characteristic zero, and hence admits for any $N>0$ a formal
logarithmic isomorphism to a power of the additive group $$ \phi : \hat{G}_{{\bf 0}}(p^N \Z_p) \cong \hat{{\bf 
G}}^d_a(p^N {\bf Z}_p)$$ preserving valuations. Hence for any $n$, 
\begin{eqnarray*} 
v(t_i(nP)) &=& v(\phi(t_i(nP))) = v(n(\phi(t_i(P))) = v(n) + v(\phi(t_i(P)) \\ &=& 
v(n) + v(t_i(P)),
\end{eqnarray*}
and we can find $n$ such that $v(t_i(P))$ becomes arbitrary large as in 
\ref{ec-3}. \qed

\vv

(c) Since by assumption, 
all elements of $A_G$ are in $\Z$, we can set $\til{a}=a$, $\ell=n$ for any $a 
\in A_G$ to get the required norm-boundedness. 

\vv

{\bf Remarks.} \ (i) From available computer algebra, the construction of 
elliptic curves which fit the above can be automated. One can compute ranks of 
elliptic curves
over $\Q$ quite fast using {\tt mwrank} by J.\ Cremona \cite{Cremona}, and over 
number fields using the {\tt gp}-package of D.\ Simon \cite{Simon}. One finds 
for example unconditionally that 
the curve $y^2=x^3+8x$ has rank one over $\Q$ and this rank stays the same
over $\Q(\sqrt[3]{2}), \Q(\sqrt[4]{2})$, hence the diophantine problem for the integers 
in these number fields has a negative answer (note that their Galois closures are non-abelian). 

\vv

(ii) We ask: given $K$, can one construct in some clever way a curve $C$ over $\Q$ such 
that its Jacobian satisfies the above conditions?

\vv

\paragraph\label{da2} {\bf Rank-preservation over $\Z$.} A similar construction (of which we leave out the details) can be performed if 
there exists a commutative (not necessarily complete)
group variety $G$ over $\Z$ such that $G(\O)$ is finitely generated and such that
$ \rk G(\Z) = \rk G(\O) > 0.$ As an example of this, let $L$ be another number 
field, linearly disjoint from $K$. Let $\langle a_i \rangle$ denote a $\Z$-basis 
for $L/\Q$ (this is also a basis for $\O_{KL}$ over $\O_K$). Let $T_L$ denote 
the
norm one torus $N^L_\Q(\sum a_i x_i) = 1$. Then $T_L(\Z) \cong \O_L^*$ and
$$ T_L(\O_K) = \ker (N_K^{KL} : \O_{KL}^* \rightarrow \O^*_K), $$
hence (by surjectivity of the relative norm) $\rk T_L(\O_K) = \rk \O^*_{KL} - 
\rk \O^*_K$. In particular, $T_L(\O_K) = T_L(\Z)$ iff
$$ r_{KL} + s_{KL} = r_K + s_K + r_L + s_L -1 $$
where $r_M,s_M$ denote the number of real, respectively half the number of complex embeddings 
of a number field $M$. 

(a) If $K$ is totally real of degree $d$, $r_{KL}=dr_L, s_{KL}=ds_L$, hence we 
want
$r_L+s_L=1$, which we can achieve by choosing $L$ quadratic imaginary; but then $T_L(\Z)$ is of rank zero.

(b) If $K$ is totally complex of degree $d$, $r_K=0,s_K=\frac{d}{2}$. Also $KL$ 
is then totally complex, so $r_{KL}=0, s_{KL} = \frac{d}{2}[L:\Q]$, hence we
want $\frac{d}{2} ([L:\Q]-1)  = r_L + s_L -1,$ but since the right hand side
is less than or equal to $[L:\Q]-1$, we find $d \leq 2$. Hence this strategy 
only works for $K$ quadratic imaginary. 

The conclusion is that this approach covers Denef's result 
 from \cite{Denef}, except that he can discard the first 
condition in our theorem 
(existence of elliptic curve of rank one) by using a torus instead. 

\vv

{\bf Remark.} In all these examples, division-ample sets are actually subsets of 
the integers. Can one find a division-ample set which does not consists 
of just ordinary integers? 

\vv

\vv

\begin{small}

\noindent {\bf Aknowledgements.} \ The authors thank Jan Van Geel for very useful help and encouragement. The third author was supported by a 
Marie-Curie Individual Fellowship  (HPMF-CT-2001-01384). 

\vv

\def\refname{\noindent\normalsize{References}}

\vv

\vv

\noindent Mathematisch Instituut,
Universiteit Utrecht, Postbus 80010, 3508 TA  Utrecht, Nederland 
(gc)\footnote{All correspondence should be sent to this author}

\vv

\noindent Department of Mathematics, University of Crete, P.O.\ Box 1470, 
Herakleio, Crete, Greece (tp)

\vv

\noindent Equipe de Logique Math\'ematique, U.F.R. de Math\'ematiques (case 
7012), Universit\'e Denis-Diderot Paris 7, 2 place Jussieu, 75251 Paris Cedex 
05, France (kz)

\vv

\noindent Email: 

{\tt cornelissen@math.uu.nl}, 

{\tt pheidas@math.uoc.gr}, 

{\tt 
zahidi@logique.jussieu.fr}

\end{small}

\end{document}